\documentclass{article}

\def \N {\mathbf N}

\def \R {{\bf R}}
\def \P {{\bf P}}
\def \T {{\bf T}}
\def \s {{\bf S}}

\def\uu{\bigsqcup}
\def\eps{\varepsilon}

\title{ \bf Non-recurrence and divergent  $\bf ({p(n)},{q(n)})$-averages \\ for deterministic automorphisms}

\author{ V.V. Ryzhikov}
\date{}

\begin{document}
\Large

\maketitle

\begin{abstract}

If sequences ${ p(n)},{ q(n)}$ satisfy
the condition $ p(n+1)- p(n), q(n+1)- q(n)\ \to\ +\infty,$  there exist automorphisms $S,T$ with simple continuous spectrum such that
the sequence
$$ \frac 1 N \sum_{n=1}^N \int S^{ p(n)}f\cdot T^{ q(n)}g \, d\mu$$
diverges for some functions $f,g\in L_2$.
In particular, this gives an answer to the question of Frantzikinakis and Host for  the case $p(n)=n^2, q(n)=n^3$. 
The automorphisms $S,T$ are realized as Gaussian and Poisson suspensions over the automorphisms $\s,\T$ of a space with sigma-finite measure. The automorphism $\T$ lies in the full group $[\s]$ of the automorphism $\s$. We give also  sketch-example of linear non-recurrence for a pair of mixing suspensions  of zero entropy and with singular and Lebesgue parts in their spectra.

\vspace{2mm}
\end{abstract}

\section{Introduction}
In  \cite{F} it  was asked about the convergence in the space $ L_2(X, \mu)$ of averages for a sequence of the form
$\{S^nf\, T^ng\}$, where $f,g\in L_\infty(\mu)$, where $S,T$ are deterministic automorphisms of the probability space $(X,\mu)$. The authors of \cite{F} proved that for such $S,T$ with $p\geq 2$ the convergence of averages of the form
$$ \frac 1 N \sum_{n=1}^N S^nf \cdot T^{n^p}g.$$ holds.
In the papers \cite{A}, \cite{Ye}, \cite{R} suitable pairs $T,S$ and $f,g\in L_\infty(X, \mu)$ with a divergent sequence $$ \frac 1 N \sum_{n=1}^N \int T^nf  \cdot S^ng \, d\mu$$ were constructed in different ways.

In \cite{F}   the question was asked also about the convergence of averages for the sequence $S^{n^2}f\, T^{n^3}g$. See   \cite{Ye24},\cite{K} for some related results. In connection with this, we  consider a more general formulation of this problem: 

\it
for which automorphisms $S,T$ with zero entropy and sequences ${ p (n)},{ q (n)}\to\infty$ can the averages $ \sum_{n=1}^N \int S^{ p(n)}f \cdot  T^{ q(n)}g \, d\mu/N$ diverge?\rm

\vspace{7mm}
\bf Theorem 1. \it
Let  sequences ${ p(n)},{ q(n)}$ satisfy
the property  $$ p(n+1)- p(n), \ \ q(n+1)- q(n)\ \to\ +\infty. $$
Then there exist automorphisms $S,T$ with simple singular spectrum and a set $C$ such that the sequence $ \sum_{n=1}^{N} \mu(S^{ p(n)}C\cap T^{ q(n)}C)/N$
diverges. \rm

\section{ Construction of  infinite automorphism $\s$}
To prove the theorem, we will use constructions of rank one, which are useful for various problems. With their help the following were solved (see, for example, \cite{R21}, \cite{R22}, \cite{R23}):

-- Kolmogorov's problem on the group property of  spectrum and Rokhlin's problem on the homogeneous spectrum of multiplicity 2 in the class of mixing automorphisms;

- ergodic homoclinic groups of an automorphism with  simple singular spectrum were discovered (Gordin's problem);

-- rigid transformations with a given mixing sequence of zero density were found (Bergelson's problem);

-- measures on a line were constructed whose tensor square has an absolutely continuous projection (onto a fixed line in the plane) along a dense set of directions and a singular projection along another dense set of directions (Oseledets's problem);

-- all sets of spectral multiplicities of the form $M\cup \{\infty\}$ ($M\subset \N$) were realized.

--  To show slow decay of correlations
for the generic  mixing automorphisms in Alperm-Tikhonov space we  used rank one constructions too. ( Podvigin's problem, see \cite{mz24}). 

-- To answer the question of Frantzikinakis and Host \cite{F}, we construct similarly to \cite{R} infinite transformations $\s, \T$ for which the corresponding Poisson (and Gaussian) suspensions will play the role $S,T$ from Theorem 1.

\bf Rank-one constructions. \rm 
Let the parameters $h_1=1$, $r_j\geq 2$ and spacer sequence
$$ \bar s_j=(s_j(1), s_j(2),\dots,s_j(r_j)), \ s_j(i)\geq 0, \ j\in\N$$ are given.
The phase space $X$ for the transformation $ \s$ defined below is the union of towers $$X_j=\uu_{i=0}^{h_j-1} \s^iE_j,$$
where $ \s^iE_j$ are disjoint half-intervals called floors.

At stage $j$, the transformation $ \s$ is defined as a usual interval transfer, but at the top floor  $ \s^{h_j-1}E_j$ of tower $X_j$, the transformation is not yet defined.
Tower $X_j$ is cut into $r_j$ identical narrow subtowers $X_{j,i}$ (they are called columns), and $s_j(i)$ new floors are added above each column $X_{j,i}$. The transformation $ \s$ is still defined as moving up one floor, but the last floor added above column $i$, the transformation $ \s$ sends to the bottom floor of column $i+1$. Thus, a new tower arises
$$X_{j+1}=\uu_{i=0}^{h_{j+1}-1} \s^iE_{j+1}$$
of height
$$h_{j+1}=r_jh_j+\sum_{i=1}^{r_j}s_j(i),$$
where $E_{j+1}$ is the lower floor of column $X_{j,1}$.
Note that by defining the transformation, we completely preserve the previous constructions.
Continuing this process to infinity,
we obtain an invertible transformation $ \s:X\to X$, preserving the Lebesgue measure on the union $X=\uu_j X_j$.
Such a transformation is called an automorphism of rank one. It has  simple
 spectrum.
\newpage
\bf Parameters of   $\s$.
\rm Let us consider a construction $\s$ with parameters
$r_j=2$, $ s_j(1)=0$, $ s_j(2)\to \infty$  (similar  constructions have been used for other purposes in \cite{R21}, \cite{R22},\cite{R23}.)  Given  $\{p(n)\},\{q(n)\}$  choose 
$s_j(2)>\max\{p(jh_j),q(jh_j)\}$.   

Let the sets $A,B$ consist of the floors of the tower of stage $j_0$,
then for $j>j_0$ it holds (this is easy to check)
$$\mu(\s^{h_j}A\cap B)=\mu(A\cap B)/2.$$
Therefore, $ \s^{h_j}\to_w I/2$, which means that  spectrum of the operator
$$\exp(\s)= {\bf 1} \oplus \s\oplus\s^{\odot 2}\oplus\s^{\odot 3}\oplus\dots \ $$
singular. Therefore, the Gaussian and Poisson suspensions over $\s$ (as operators isomorphic to the operator $\exp(\s)$) have singular spectrum. Their entropy is zero. Later we will say how to modify the construction of $\s$ to provide  simple spectrum of suspensions.

\vspace{3mm}

\section{ Proof of the theorem}
\bf Definition of  $\T$. \rm The automorphism $\s$ (almost) cyclically permutes the  floors in the towers $X_j$.  
Our $\T$ will   permutate these floors in a different order. We define $\T$ as $\P\s\P^{-1}$, where $\P$   permutes the floorsas well. Let us describe it.

We identify the set of floors of the tower of stage $j+1$ with the ordered
set of numbers $I_j=(0,1,2,\dots, h_{j+1}-1)$. The automorphism $\s$
obviously corresponds to an (almost) cyclic permutation $\sigma$ on $I_j$. It is given by the formula $\sigma(n)=n+1$, but let the value $\sigma(h_{j+1}-1)$ be undefined.
(In the spirit of Katok-Stepin approximations, we can formally set $\sigma(h_{j+1}-1)=0$.)

For $n>n_j$, the intervals $(p(n),p(n)+1,p(n)+2, \dots, p(n)+2h_j-1)$ 
do not intersect as well as   the intervals
$(q(n)+2h_j,q(n)+2h_j+1,q(n)+2h_j+2, \dots,q(n) +4h_j-1)$.

\newpage
\bf For odd $\bf j$ \rm we find a permutation $\pi_j$ on $\{2h_j, 2h_j-1,
\dots, h_{j+1}-1\}$   such that
$$\pi_j(q(n),q(n)+1,q(n)+2, \dots, q(n)+2h_j-1)=
(p(n)+2h_j,q(n)+2h_j+1, \dots,q(n) +4h_j-1),$$
$ p(n),q(n)>2h_j$, $p(n),q(n)< h_{j+1} - 4h_j$. 
With this permutation $\pi_j$ we associate an automorpism $\P_j$, which exchanges floors, sending   $k$-floor  by a parallel translation to  $\pi_j(k)$-floor.
Outside of $X_{j+1}\setminus X_j$  our $\P_j$
acts identically. Let's recall that  $$X_{j}=\uu_{i=0}^{2h_{j}-1} \s^iE_{j+1}, \ \
X_{j+1}=\uu_{i=0}^{h_{j+1}-1} \s^iE_{j+1}.$$

Let $A=X_1\subset X_j$, then $A$ consists of the floors of stage $j+1$.
(we  identify $A$ with the set of numbers of these floors, which form a subset of  $(0,1,2, \dots, 2h_j)$.)
Put $\T=\P_j\s\P_j^{-1}$, then
$$\T^{q(n)}A= \s^{p(n)+2h_j}A,$$
$$\T^{q(n)}A\cap \s^{p(n)}A=\emptyset, \eqno (1) $$

\bf For even $\bf j$ \rm we define a permutation $\pi_j$:
$$\pi_j(q(n),q(n)+1,q(n)+2, \dots, q(n)+2h_j-1)=
(p(n),p(n)+1, \dots,p(n) +2h_j-1),$$
$ p(n),q(n)>2h_j$, $p(n),q(n)< h_{j+1} - 2h_j$.
 The corresponding automorphism $\P_j$ again send $k$-floor
to  $\pi_j(k)$-floor, and outside of $X_{j+1}\setminus X_j$  we set  $\P_j=Id$.

For  $\T=\P_j\s\P_j^{-1}$ we have  
$$\T^{q(n)}A=\s^{p(n)}A.\eqno (2)$$

Let $$\P=\prod_j \P_j, \ \ \T=\P\s\P^{-1}.$$  It is important that equalities (1), (2)  are preserved for such $\T$.

So we can find  $N_j\to\infty$ (with corresponding $h_j\to\infty$)  such that 
$$ \frac 1 {N_{2k}} \sum_{n=1}^{N_{2k}} \mu(\s^{ p(n)}A\cap \T^{q(n)}A)\to \mu(A), $$
and 
$$ \frac 1 {N_{2k+1}} \sum_{n=1}^{N_{2k+1}} \mu(\s^{ p(n)}A\cap \T^{ q(n)}A)\to 0.$$
We choose $N_j=\max \{n\,:\,  p(n),q(n)<h_{j+1}-4h_j$.
Note that $N_{j-1}/N_{j}\to 0$ because $N_{j-1}<h_j$, $N_{j}>jh_j$  (recall that we choose    above $h_{j+1}>\max\{p(jh_j),q(jh_j)\}$). 

\bf Modification. \rm For some rare set of  stages $j$, we can modify the parameters of the above construction $\s$ so that the weak closure of its powers will possesses all possible polynomial. For this let $r_j=2j$,  $s_j(i)=0$ as $i\leq j$, and  $s_j(i)=1$ as $i\leq 2j$
(Katok's spacers).    This   provides  simple spectrum  for the operator $\exp(\s)$. At such stages we set $\P_j=Id$.

\bf  Completion of the proof. \rm Let us consider  in the Poisson space  the cylinder $C=C(A,0)$ (see \cite{R}),  $0<c=\mu_\circ(C)<1$,  $\mu_\circ$ is the Poisson measure. If sets $A,A'$ do not intersect, then the corresponding cylinders $C,C'$ are independent
with respect to $\mu_\circ$.
Since $S^n C(A,0)=C(\s^nA,0)$,
$T^n C(A,0)=C(\T^nA,0)$, we obtain
$$ \frac 1 {N_{2k+1}} \sum_{n=1}^{N_{2k+1}} \mu_\circ(S^{p(n)}C\cap T^{q(n)}A)\to c^2,$$
$$ \frac 1 {N_{2k}} \sum_{n=1}^{N_{2k}} \mu_\circ(S^{p(n)}C\cap T^{q(n)}A)\to c.$$

\section{Linear non-recurrence}
In \cite{F} was stated  a problem on recurrence. \it
Are there deterministic automorphisms $S,T$ such that for a set $D$ of positive measure $S^nD\cap T^nD=\emptyset$ for all $n>0$? \rm

 We have already mentioned Gordin's problem on homoclinic groups.  Let $T$ be an automorphism of a probability space, the Gordin homoclinic group $H(T)$ is defined as
$$H(T)=\{R\in Aut\,:\, T^{-n}RT^n\, \to_s\, I, \ n\to\infty\}.$$
It is easy to see that  mixing    Poisson suspension over a "mixing" $\T$ ($T^n\to_w 0$) has the following homoclinic element.  
Let $\mu ( x:  \R x\neq x)<\infty$. Then $\T^{-n}\R \T^n\, \to_s\, {\bf I}$, so for the corrsponding Poisson suspensions we have
$$T^{-n}R T^n\, \to_s\, {\ I}.$$    Now let   $\T$ be a rank-one transformation, and $\R\in [\T]$ be  interval involutive exchange of some  floors $E$ and $\R E$. Outside 
$$F=E\uu \R E$$ we set $\R =Id $.  In the Poisson space, let's consider the set  $$D=C(E,0)\bigcap C(RE,1).$$ Recall that $C(E,0)$ consists of all countable sets $x_\circ=\{x_i\}$ such that $E\cap x_\circ$ has cardinality 0. (Note that $x_\circ$ has no limit points in the intervals that make up the phase space $X$.)
The cylinder $C(\R E,1)$ consists of all $x_\circ$ such that $RE\cap x_\circ$ has cardinality 1.  So $$R D\bigcap D=\emptyset.$$

Set $$S=R^{-1}TR.$$ Since $$\mu_\circ (T^{-n}RT^n D\Delta D)\to 0$$ we get
$$\mu_\circ (S^nD\cap T^n D)= \mu_\circ (RD\cap T^{-n}RT^n D)\ \to 0.$$
If the series $\sum_n \mu_\circ (RD\cap T^{-n}RT^n D)$ converges, then we achieve the goal by considering its power $T^p$ instead of $T$. For $T^p$ the sum of the corresponding series  can be arbitrarily small (if $p$ is large) and so less than the measure of the set $D$. 
Then  desired $D'\subset D$ of a positive  measure can be found:
 $$S^nD'\cap T^nD'=\emptyset, \ n>0.$$

But this is impossible.  In fact, the above  series always diverges.
However we do the following. We see that 
 $$\mu_\circ (RD\cap T^{-n}RT^n D)< const \cdot \mu(\T^n F\cap F).$$ 
Indeed, let $\mu(\T^n F\cap F)=0$ this means that $\T^{-n}\R\T^n$  is the identity  on  $F$.  Thus $T^{-n}RT^n D=D$.  If $\mu(\T^n F\cap F)=\eps>0$,
then $$\mu_\circ (RD\cap T^{-n}RT^n D)\sim const \eps.$$

The  correlations $\mu(\T^n F\cap F)$  for some Sidon constructions $\T$ can be so that  $$\sum_n \mu(\T^n F\cap F)^4<\infty$$
 (see  \cite{24msb}). Then  
$$\sum_n \mu_\circ^{\otimes 4}
((S^nD)^{\times 4}\cap (T^n D)^{\times 4})<\infty.$$
Some powers  $(S^{\times 4})^p$ and $(T^{\times 4})^p$, and 
$D'\subset D^{\times 4}$ of a positive measure $\mu_\circ^{\otimes 4}$ ensure the absence of recurrence.  
Our Poisson suspensios have  zero entropy, so get desired linear non-recurrence.

\vspace{7mm}
The author thanks Nikos Frantzikinakis for stimulating questions and pleasant discussions on the above  topic, and Jean-Paul Thouvenot for his interest and friendly support. 

\large

\end{document}